\documentclass[11pt,reqno]{amsart}
\usepackage{amssymb,amsmath}\newcommand{\be}{\begin{eqnarray}}
\newcommand{\ee}{\end{eqnarray}}

\newcommand{\Bel}{\mathbf B}

\newcommand{\sign}{\operatorname{sign}}

\newcommand{\half}{\frac12}

\newcommand{\R}{{\mathbb R}}\newcommand{\Q}{{\mathbb Q}}

\newcommand{\B}{{\mathcal B}}
\newcommand{\E}{{\bf E\,}}

\newtheorem{theorem}{Theorem}\newtheorem{lemma}[theorem]{Lemma}

\theoremstyle{definition}
\theoremstyle{remark}\numberwithin{equation}{section}\input epsf.sty
\begin{document}\thispagestyle{empty}

\title[Orthogonal martingales]{Some new Bellman functions and subordination by orthogonal  martingales in $L^{p}, 1<p\le 2$.}
\author{Prabhu Janakiraman}\address{Prabhu Janakiraman, Department of Mathematics, Michigan State University,
{\tt pjanakir1978@gmail.com}}
\author{Vasily Vasyunin}\address{Vasily Vasyunin, V.A. Steklov. Math. Inst.
{\tt vasyunin@pdmi.ras.ru}}
\author{Alexander Volberg}\address{Alexander Volberg, Department of Mathematics, Michigan State University,
{\tt volberg@math.msu.edu}}


\maketitle

\section{Introduction: Orthogonal martingales and the Beurling-Ahlfors transform}

The main result of this note is Theorem \ref{less2} below. The main interest is the array of new Bellman function, very different from Burkholder's function. 

A complex-valued martingale $Y=Y_1+iY_2$ is said to be $orthogonal$ if the quadratic variations of the coordinate martingales are equal and their mutual covariation is $0$:
\[\left<Y_1\right>=\left<Y_2\right>,\hspace{3mm} \left<Y_1,Y_2\right>=0.\]
In \cite{BaJ1}, Ba\~nuelos and Janakiraman make the observation that the martingale associated with the Beurling-Ahlfors transform is in fact an orthogonal martingale. They show that Burkholder's proof in \cite{Bu3} naturally accommodates for this property and leads to an improvement in the estimate of $\|B\|_p$.
\begin{theorem}\label{BurkBaJa} (One-sided orthogonality as allowed in Burkholder's proof)
\begin{enumerate}
\item (Left-side orthogonality) Suppose $2\leq p<\infty$. If $Y$ is an orthogonal martingale and $X$ is any martingale such that $\left<Y\right>\leq \left<X\right>$, then 
\begin{equation}\label{estBJ}
\|Y\|_p \leq \sqrt{\frac{p^2-p}{2}}\|X\|_p.
\end{equation}
\item (Right-side orthogonality) Suppose $1<p<2$. If $X$ is an orthogonal martingale and $Y$ is any martingale such that $\left<Y\right>\leq \left<X\right>$, then 
\begin{equation}\label{estBJ1}
\|Y\|_p \leq \sqrt{\frac{2}{p^2-p}}\|X\|_p.
\end{equation}
\end{enumerate}It is not known whether these estimates are the best possible.
\end{theorem}

\bigskip

\noindent{\bf Remark.}
The result for right-side orthogonality is stated in \cite{JV} and not in \cite{BaJ1}. In \cite{JV} we emulate \cite{BaJ1} to prove in a very simple way an estimate on right-side orthogonality and in the regime $1<p\le 2$. In the present work we tried to have a better constant for this regime--as the sharpness in \cite{BaJ1} and \cite{JV} is somewhat dubious. We build for that some    family of new (funny and interesting) Bellman functions very different from the original Burkholder's function.    Even though the approach is quite different from the one in \cite{BaJ1} and \cite{JV}, the constants we will obtain here are {\it the same}!
So may be they are  sharp after all.

\bigskip

If $X$ and $Y$ are the martingales associated with $f$ and $Bf$ respectively, then $Y$ is orthogonal, $\left<Y\right>\leq 4\left<X\right>$ and hence by (1), we obtain
\begin{equation}\label{Beurest}
\|Bf\|_p\leq \sqrt{2(p^2-p)}\|f\|_p \textrm{ for } p\geq 2.
\end{equation}
By interpolating this estimate $\sqrt{2(p^2-p)}$ with the known $\|B\|_2=1$, Ba\~nuelos and Janakiraman establish the present best estimate in publication:
\begin{equation}\label{best-est}
\|B\|_p\leq 1.575 (p^*-1).
\end{equation}

\section{New Questions and  Results}
Since $B$ is associated with left-side orthogonality and since we know $\|B\|_p = \|B\|_{p'}$, two important questions are
\begin{enumerate}
\item If $2\leq p<\infty$, what is the best constant $C_p$ in the left-side orthogonality problem: $\|Y\|_p\leq C_p\|X\|_p$, where $Y$ is orthogonal and $\left<Y\right>\leq \left<X\right>$?
\item Similarly, if $1<p'<2$, what is the best constant $C_{p'}$ in the left-side orthogonality problem?
\end{enumerate}
We have separated the two questions since Burkholder's proof (and his function) already gives a good answer when $p\geq 2$. It may be (although we have now some doubts about that) the best possible as well. However no estimate (better than $p-1$) follows from analyzing Burkholder's function when $1<p'<2$. Perhaps, we may hope, $C_{p'}<\sqrt{\frac{p^2-p}{2}}$ when $1<p'=\frac{p}{p-1}<2$, which would then imply a better estimate for $\|B\|_p$. This paper 'answers' this hope in the negative by finding $C_{p'}$; see Theorem \ref{MAINThm}. We also ask and answer the analogous question of right-side orthogonality when $2<p<\infty$. In the spirit of Burkholder \cite{Bu8}, we believe these questions are of independent interest in martingale theory and may have deeper connections with other areas of mathematics.

\bigskip 

\noindent{\bf Remark.}
The following sharp estimates are proved in \cite{BJV_La}, they cover the left-side orthogonality for the regime $1<p\le 2$ and the right-side orthogonality for the regime $2\le p<\infty$. Notice that two complementary regimes have the estimates: for $2\le p<\infty$ and left-side orthogonality in \cite{BaJ1}, for $1<p\le 2$ in this note and in \cite{JV}, but the sharpness is  somewhat dubious.

\bigskip

\begin{theorem}\label{MAINThm}Let $Y=(Y_1,Y_2)$ be an orthogonal martingale and $X=(X_1,X_2)$ be an arbitrary martingale.
\begin{enumerate}
\item Let $1<p' \leq 2$. Suppose $\left<Y\right>\leq\left<X\right>$. Then the least constant that always works in the inequality $\|Y\|_{p'}\leq C_{p'}\|X\|_{p'}$ is 
\begin{equation}\label{constp<2}
C_{p'} = \frac{1}{\sqrt{2}}\frac{z_{p'}}{1-z_{p'}}
\end{equation}
where $z_{p'}$ is the least positive root in $(0,1)$ of the bounded Laguerre function $L_{p'}$.
\item Let $2\leq p<\infty$. Suppose $\left<X\right>\leq \left<Y\right>$. Then the least constant that always works in the inequality $\|X\|_{p}\leq C_{p}\|Y\|_{p}$ is 
\begin{equation}\label{constp>2}
C_{p} = \sqrt{2}\frac{1-z_{p}}{z_{p}}
\end{equation}
where $z_{p}$ is the least positive root in $(0,1)$ of the bounded Laguerre function $L_{p}$.
\end{enumerate}
\end{theorem}
The Laguerre function $L_p$ solves the ODE 
\[sL_p''(s)+(1-s)L_p'(s)+pL_p(s) = 0.\] These functions are discussed further and their properties deduced in section \eqref{laguerresection}; see also \cite{BJV}, \cite{C}, \cite{CL}.

As mentioned earlier, (based however on numerical evidence) we believe in general $\sqrt{\frac{p^2-p}{2}}<C_{p'}<p-1$ and that these theorems cannot imply better estimates for $\|B\|_p$. However based again on numerical evidence, the following conjecture is made.

\noindent{\bf Conjecture.}
For $1<p'=\frac{p}{p-1}<2$, $C_{p'} = C_p$, or equivalently,
\[ \frac{1}{\sqrt{2}}\frac{z_{p'}}{1-z_{p'}} = \sqrt{2}\frac{1-z_{p}}{z_{p}}.\]

It is conjecture relating the roots of the Laguerre functions. Notice that such a statement is not true with the constants from Theorem \ref{BurkBaJa}, and $\sqrt{\frac{2}{p'^2-p'}}<\sqrt{\frac{p^2-p}{2}}$ for all $p>2$. So this conjecture (if true) suggests some distinct implications for the two settings. Note on the other hand, that the form of the two sets of constants are very analogous.

\section{Orthogonality}

Let $Z=(X,Y), W=(U, V)$ be two $\R^2$-valued martingales on the filtration of $2$--dimensional Brownian motion $B_s= (B_{1s}, B_{2s})$. Let $A=\begin{bmatrix} -1, & i\\
i, & 1\end{bmatrix}$. We want $W$ to be a martingale transform of $Z$ defined by $A$. Let
 $$X(t) = \int_0^t \overrightarrow{x}(s)\cdot dB_s\,,$$ 
$$Y(t) = \int_0^t\overrightarrow{y}(s)\cdot dB_s\,,$$
 where $X, Y$ are {\it real-valued} processes, and $\overrightarrow{x}(s), \overrightarrow{y}(s)$ are $\R^2$--valued ``martingale differences". 
 
 Put
 \begin{equation}
\label{G}
Z(t) = X(t) +iY(t)\,, Z(t) = \int_0^t (\overrightarrow{x}(s) +i \overrightarrow{y}(s))\cdot dB_s\,,
\end{equation}
and
\begin{equation}
\label{F}
W(t) = U(t) +iV(t)\,, W(t) = \int_0^t (A(\overrightarrow{x}(s) +i \overrightarrow{y}(s)))\cdot dB_s\,.
\end{equation}
We will denote
$$
W=A\star Z\,.
$$
As before 
$$
U(t) = \int_0^t \overrightarrow{u}(s)\cdot dB_s\,,
$$ 
$$
V(t) = \int_0^t\overrightarrow{v}(s)\cdot dB_s\,,
$$
$$
W(t) = \int_0^t (\overrightarrow{u}(s) +i \overrightarrow{v}(s))\cdot dB_s\,.
$$

 We can easily write components of $\overrightarrow{u}(s), \overrightarrow{v}(s)$:
$$
u_1(s) =-x_1(s)-y_2(s)\,,\,\,v_1(s) = x_2(s)-y_1(s)\,, i=1,2\,,
$$
$$
u_2(s) =x_2(s)-y_1(s)\,,\,\,v_2(s) =x_1(s)+y_2(s)\,, i=1,2\,.
$$

Notice that
\begin{equation}
\label{orto}
\overrightarrow{u}\cdot \overrightarrow{v} = u_1v_1+u_2v_2 = -(x_1+y_2)(x_2-y_1) + (x_2-y_1)(x_1+y_2)=0\,.
\end{equation}

\subsection{Local ortogonality.}
\label{localort}

The processes

$$
\langle X, U\rangle(t) :=\int_0^t \overrightarrow{x}\cdot \overrightarrow{u}ds\,,
\,\,
\langle X, V\rangle(t) :=\int_0^t \overrightarrow{x}\cdot \overrightarrow{v}ds\,,
$$
$$
\langle Y, U\rangle(t) :=\int_0^t \overrightarrow{y}\cdot \overrightarrow{u}ds\,,
\,\,
\langle Y, V\rangle(t) :=\int_0^t \overrightarrow{y}\cdot \overrightarrow{v}ds\,,
$$
$$
\langle X, X\rangle(t) :=\int_0^t \overrightarrow{x}\cdot \overrightarrow{x}ds\,,
\,\,
\langle Y, Y\rangle(t) :=\int_0^t \overrightarrow{y}\cdot \overrightarrow{y}ds\,,
$$
$$
\langle X, Y\rangle(t) :=\int_0^t \overrightarrow{x}\cdot \overrightarrow{y}ds\,,
\,\,
\langle U, U\rangle(t) :=\int_0^t \overrightarrow{u}\cdot \overrightarrow{u}ds\,,
$$
$$
\langle V, V\rangle(t) :=\int_0^t \overrightarrow{v}\cdot \overrightarrow{v}ds\,,
\,\,
\langle U, V\rangle(t) :=\int_0^t \overrightarrow{u}\cdot \overrightarrow{v}ds\,.
$$
are called the covariance processes. We can denote

$$
d\langle X, U\rangle(t) := \overrightarrow{x}(t)\cdot \overrightarrow{u}(t)\,,
\,\,
d\langle X, V\rangle(t) := \overrightarrow{x}(t)\cdot \overrightarrow{v}(t)\,,
$$
$$
d\langle Y, U\rangle(t) :=\overrightarrow{y}(t)\cdot \overrightarrow{u}(t)\,,
\,\,
d\langle Y, V\rangle(t) :=\overrightarrow{y}(t)\cdot \overrightarrow{v}(t)\,,
$$
$$
d\langle X, X\rangle(t) :=\overrightarrow{x}(t)\cdot \overrightarrow{x}(t)\,,
\,\,
d\langle Y, Y\rangle(t) := \overrightarrow{y}(t)\cdot \overrightarrow{y}(t)\,,
$$
$$
d\langle X, Y\rangle(t) := \overrightarrow{x}(t)\cdot \overrightarrow{y}(t)\,,
\,\,
d\langle U, U\rangle(t) :=\overrightarrow{u}(t)\cdot \overrightarrow{u}(t)\,,
$$
$$
d\langle V, V\rangle(t) :=\overrightarrow{v}(t)\cdot \overrightarrow{v}(t)\,,
\,\,
d\langle U, V\rangle(t) :=\overrightarrow{u}(t)\cdot \overrightarrow{v}(t)\,,
$$
$$
d\langle Z, Z\rangle(t) := (\overrightarrow{x}(t)\cdot \overrightarrow{x}(t)+\overrightarrow{y}(t)\cdot \overrightarrow{y}(t))\,,
\,\,
d\langle W, W\rangle(t) := (\overrightarrow{u}(t)\cdot \overrightarrow{u}(t)+\overrightarrow{v}(t)\cdot \overrightarrow{v}(t))\,.
$$

Important is an observation

\begin{lemma}
\label{locort}
Let $A=\begin{bmatrix} -1, & i\\
i, & 1\end{bmatrix}$. Then
\begin{equation}
\label{ortheq1}
d\langle U, V\rangle (t) =0\,.
\end{equation}
Or
$$
\overrightarrow{u}(t)\cdot \overrightarrow{v}(t)=0\,.
$$ 
\end{lemma}

Also

\begin{lemma}
\label{subord4}
With the same $A$
$$
d\langle U, U\rangle (t) \le 2\,d\langle Z, Z\rangle (t)\,.
$$
$$
d\langle V, V\rangle (t) \le 2\,d\langle Z, Z\rangle (t)\,.
$$
Or 
$$
\overrightarrow{u}(t)\cdot \overrightarrow{u}(t)\le 2\,(\overrightarrow{x}(t)\cdot \overrightarrow{x}(t) + \overrightarrow{y}(t)\cdot \overrightarrow{y}(t))\,,
$$ 
$$
\overrightarrow{v}(t)\cdot \overrightarrow{v}(t)\le 2\,(\overrightarrow{x}(t)\cdot \overrightarrow{x}(t) + \overrightarrow{y}(t)\cdot \overrightarrow{y}(t))\,.
$$ 
Or
\begin{equation}
\label{subord4eq}
d\langle W, W\rangle (t) \le 4\,d\langle Z, Z\rangle (t)\,.
\end{equation}
\end{lemma}

\begin{proof}
$$
\overrightarrow{u}(t)\cdot \overrightarrow{u}(t)= (x_1 + y_2)^2 + (x_2-y_1)^2=2\,(x_1\,y_2 - x_2\, y_1) +
$$
$$
(x_1)^2 + (y_2)^2 + (x_2)^2 + (y_1)^2 \le 2\,((x_1)^2 + (y_2)^2 + (x_2)^2 + (y_1)^2 )= 2\,d\langle Z, Z\rangle\,.
$$
The same for $v$.

\end{proof}

\medskip

\noindent{\bf Definition.} The complex martingale $W=A\star Z$ will be called the Ahlfors-Beurling transform of martingale $Z$.

Now let us quote a theorem of Banuelos-Janakiraman \cite{BaJ1}:

\begin{theorem}
\label{BaJthMart}
Let $Z,W$ be two martingales on the filtration of Brownian motion, let $W$ be an orthogonal martinagle in the sense of \eqref{ortheq1}: $d\langle U, V\rangle =0$ and there is a subordination property
\begin{equation}
\label{subord1eq2}
d\langle W, W\rangle \le d\langle Z,Z\rangle
\end{equation}
Let $p\ge 2$. Then for every $t$ ($|\cdot|$ denotes the euclidean norm in $\R^2$)
\begin{equation}
\label{pnorm1}
(\E |W|^p)^{1/p} \le \sqrt{\frac{p^2-p}{2}}(\E |Z|^p)^{1/p}\,.
\end{equation}
\end{theorem}

We will use the notations
$$
\|Z\|_p:= (\E |Z|^p)^{1/p}\,.
$$

Applied to our case (with the help of Lemmas \ref{locort}, \ref{subord4}) we get from Theorem \ref{BaJthMart} the following

\begin{theorem}
\label{FG1}
$\|W\|_p =\|A\star Z\|_p \le \sqrt{2(p^2-p)} \|Z\|_p\,,\,\,\forall p\ge 2\,.$
\end{theorem}

\section{Subordination by orthogonal martinagales $L^{3/2}$}
\label{L32}


For $1<p\le 2$ one has the following

\begin{theorem}
\label{less2}
Let $Z, W$ be two $\R^2$ martingales as above, and $W$ is an orthogonal martingale 
:
$$
d\langle U, V\rangle =0\,,
$$
satisfying also
\begin{equation}
\label{equalnorms1}
d\langle U, U\rangle =d\langle V, V\rangle\,.
\end{equation}
Let $Z$ be subordinated to the orthogonal martingale $W$:
\begin{equation}
\label{subord1}
d\langle Z, Z\rangle\le \langle W, W\rangle 
\end{equation}
Then for $1<q\le 2$
\begin{equation}
\label{less2eq}
\|Z\|_q\le \sqrt{\frac{2}{q^2-q}}\|W\|_q\,.
\end{equation}
\end{theorem}

We will give a proof, but first it will be given  for $q=3/2$, and only later for all $q\in (1, 2]$. Moreover our proof may indicate--especially compared with a completely different proof having the same result in \cite{JV}-- that the constant $\sqrt{\frac{2}{p^2-p}}$ is sharp after all. But we cannot be sure.

\begin{proof}
We can assume that $F=(\Phi, \Psi)$ (or $F=\Phi+i\Psi$) is a martingale on the filtration of Brownian motion
$$
\Phi(t) = \int_0^t \overrightarrow{\phi}(s)\cdot dB_s\,,
$$ 
$$
\Psi(t) = \int_0^t\overrightarrow{\psi}(s)\cdot dB_s\,,
$$$$
X(t) = \int_0^t \overrightarrow{x}(s)\cdot dB_s\,,
$$ 
$$
Y(t) = \int_0^t\overrightarrow{y}(s)\cdot dB_s\,,
$$$$
U(t) = \int_0^t \overrightarrow{u}(s)\cdot dB_s\,,
$$ 
$$
V(t) = \int_0^t\overrightarrow{v}(s)\cdot dB_s\,,
$$
and that these vector processes and their components satisfy
Lemmas \ref{locort} and \ref{subord4}:
\begin{equation}
\label{locort1}
u_1 v_1 + u_2v_2 = 0\,,
\end{equation}
\begin{equation}
\label{equalnorms2}
(u_1)^2 + (u_2)^2=(v_1)^2 + (v_2)^2\,,
\end{equation}

$$
\Im \E (F \cdot Z )= \int_0^t (d\langle \Phi, X\rangle + d\langle \Psi, Y\rangle)\,ds =
\int_0^t (\phi_1 x_1 + \phi_2 x_2 + \psi_1 y_1 +\psi_2 y_2) \,ds\,.
$$
Hence
\begin{equation}
\label{sq0}
|\Im \E (Z \cdot F)|\le \int_0^t ((\phi_1)^2 + (\phi_2)^2 + (\psi_1)^2 +(\psi_2)^2)^{1/2}((x_1)^2 + (x_2)^2 + (y_1)^2 +(y_2)^2)^{1/2} \,ds\,.
\end{equation}

By subordination assumption \eqref{subord1} we have

\begin{equation}
\label{sq}
|\Im \E (Z \cdot F)|\le \int_0^t ((u_1)^2 + (u_2)^2 + (v_1)^2 +(v_2)^2)^{1/2}((\phi_1)^2 + (\phi_2)^2 + (\psi_1)^2 +(\psi_2)^2)^{1/2} \,ds\,.
\end{equation}

 Our next goal will be to prove that
 $$
 \sqrt{\frac32}\int_0^t((u_1)^2 + (u_2)^2 + (v_1)^2 +(v_2)^2)^{1/2}((\phi_1)^2 + (\phi_2)^2 + (\psi_1)^2 +(\psi_2)^2)^{1/2} \le
 $$
 \begin{equation}
 \label{sqrt32}
  2\bigg(\frac{\|W\|_{3/2}^{3/2}}{3/2} + \frac{\|F\|_3^3}{3}\bigg)\,.
 \end{equation}
 
 Polarize the last equation to convert its RHS to $2\|W\|_{3/2}\|F\|_3$. Then use the combination of \eqref{sq} and \eqref{sqrt32}.  Then we obtain the desired estimate
 \begin{equation}
 \label{32}
 \|Z\|_{3/2} \le \frac{2\sqrt{2}}{\sqrt{3}}\|W\|_{3/2}\,,
 \end{equation} 
 which we saw is equivalent to the claim of Theorem \ref{less2} for $q=3/2$.
 
 We are left to prove \eqref{sqrt32}. For that we will need the next section.

\end{proof}

\section{Bellman functions and Martinagales}
\label{bfm1}

Suppose we have the function of $4$ real variables such that

\begin{equation}
\label{pr1}
B(y_{11}, y_{12}, y_{21}, y_{22}) \le \frac{2}{3} (y_{11}^2 + y_{12}^2)^{3/2} + \frac{4}{3}(y_{21}^2+y_{22}^2)^{1/2}\,,
\end{equation}

\begin{equation}
\label{pr2}
\langle d^2 B(y_{11}, y_{12}, y_{21}, y_{22})\begin{bmatrix} dy_{11}\\dy_{12}\\dy_{21}\\dy_{22}\end{bmatrix},\begin{bmatrix} dy_{11}\\dy_{12}\\dy_{21}\\dy_{22}\end{bmatrix}\rangle\ge
\end{equation}
\begin{align*}
\tau(dy_{11}^2+dy_{12}^2) + \frac1\tau (dy_{21}^2+dy_{22}^2)
+\frac{3\tau}{4x_2} \Big(\frac{y_{22}dy_{21}-y_{21}dy_{22}}{x_2}\Big)^2
\\
+\frac{\tau x_1}{\sqrt{x_1^2+3x_2}}
\Big[\frac{y_{11}dy_{11}+y_{12}dy_{12}}{x_1}+
\frac1\tau\frac{y_{21}dy_{21}+y_{22}dy_{22}}{x_2}\Big]^2\,,
\end{align*}
where
\begin{equation}
\label{tau}
\frac{3}{4}\frac{\tau}{(y_{21}^2 +y_{22}^2)^{1/2}} +\frac2{\tau} \ge  \frac3{\tau}\,.
\end{equation}

Then we can prove \eqref{sqrt32}. Let us start by writing It\^o's formula for the process $b(t):=B(\Phi(t), \Psi(t), U(t),V(t))$:
$$
db = \langle \nabla B (\Phi,..,V), (d\Phi(t),..., dV(t)\rangle + \frac12(d^2 B(\phi_1,\psi_1,u_1, v_1) +d^2 B(\phi_2,\psi_2,u_2, v_2) )\,.
$$

Here $d^2B $ stands for the Hessian bilinear form. It is applied to vector $(\phi_1,\psi_1,u_1, v_1)$ and then to vector $(\phi_2,\psi_2,u_2, v_2)$. Of course second derivatives of $B$ consituting this form are calculated at point $(\Phi,\Psi,U,V)$. All this is at time $t$. The first term is a martingale with zero average, and it disappears after taking the expectation.

Therefore,
$$
\E(b(t)-b(0)) =\E\int_0^t db(s) \,ds =
$$
\begin{equation}
\label{B1}
\frac12\int_0^t ( (d^2 B(\phi_1,\psi_1,u_1, v_1) +d^2 B(\phi_2,\psi_2,u_2, v_2) )\,ds=: \frac12\int_0^t dI\,.
\end{equation}
The sum in \eqref{B1} is the Hessian bilinear form on vector $(\phi_1, \psi_1, u_1, v_1)$ plus the Hessian bilinear form on vector $(\phi_2, \psi_2, u_2, v_2)$.  Using \eqref{pr2} we can sum up these two forms with a definite cancellation:
$$
dI = \tau ((\phi_1)^2+( \psi_1)^2) +1/\tau (( u_1)^2+( v_1)^2) + \frac{3}{4} \frac{\tau}{(U^2 +V^2)^{1/2}}
\frac{V^2 (u_1)^2 + U^2 (v_1)^2 -2 UVu_1v_1}{U^2 +V^2} +\,\text{Posit.} 
$$
$$
+\tau ((\phi_2)^2+( \psi_2)^2) +1/\tau (( u_2)^2+( v_2)^2) + \frac{3}{4} \frac{\tau}{(U^2 +V^2)^{1/2}}
\frac{V^2 (u_2)^2 + U^2 (v_2)^2 -2 UVu_2v_2}{U^2 +V^2} +\,\text{Posit.}\,.
$$
 
Notice that orthogonality \eqref{locort1} and equality of norms \eqref{equalnorms2}:
\begin{equation}
\label{locort2}
d\langle U, V\rangle = 0\,,
\end{equation}
\begin{equation}
\label{equalnorms3}
d\langle U, U\rangle = d\langle V, V\rangle \,,
\end{equation}
imply pointwise equalities $u_1v_1+u_2v_2=0$ and 
$$
V^2 (u_1)^2 + U^2 (v_1)^2 +V^2 (u_2)^2 + U^2 (v_2)^2=\frac12(U^2+V^2)((u_1)^2 + (u_2)^2 +(v_1)^2 + (v_2)^2)\,.
$$
Therefore, $UV$--term will disappear, and we will get
$$
dI = \tau ((\phi_1)^2+( \psi_1)^2+(\phi_2)^2+( \psi_2)^2) +1/\tau (( u_1)^2+( v_1)^2+( u_2)^2+( v_2)^2)) + 
$$
$$
 \frac{3}{4} \frac{\tau}{(U^2 +V^2)^{3/2}}\cdot \frac12(U^2+V^2)((u_1)^2 + (u_2)^2 +(v_1)^2 + (v_2)^2)+\,\text{Posit.}=
 $$
 $$
\tau ((\phi_1)^2+( \psi_1)^2+(\phi_2)^2+( \psi_2)^2) + \frac12 \bigg(\frac{3}{4} \frac{\tau}{(U^2 +V^2)^{1/2}}+\frac2{\tau}\bigg) (( u_1)^2+( v_1)^2+( u_2)^2+( v_2)^2)) +\,\text{Posit.}
\,.
$$

Hence, by using \eqref{tau} we get
$$
dI \ge \tau (\|\overrightarrow{\phi}\|^2 +\|\overrightarrow{\psi}\|^2) + \frac32\frac1{\tau}(\|\overrightarrow{u}\|^2+\|\overrightarrow{v}\|^2)\,. 
$$
\begin{equation}
\label{dI}
\ge 2\sqrt{\frac32} (\|\overrightarrow{\phi}\|^2 +\|\overrightarrow{\psi}\|^2)^{1/2}(\|\overrightarrow{u}\|^2+\|\overrightarrow{v}\|^2)^{1/2}\,.
\end{equation}

Let us combine now \eqref{B1} and \eqref{dI}. We get

\begin{equation}
\label{B1dI}
\sqrt{\frac32} \int_0^t(\|\overrightarrow{\phi}\|^2 +\|\overrightarrow{\psi}\|^2)^{1/2}(\|\overrightarrow{u}\|^2+\|\overrightarrow{v}\|^2)^{1/2}\,ds \le \frac{1}{2} dI \le \E (b(t))\,.
\end{equation}
We used  \eqref{pr1} that claims  $b\ge 0$. But it also claims that
\begin{equation}
\label{bsverhu}
b(t) = B(\Phi(t), \Psi(t),U(t), V(t))\le 2\,\bigg(\frac{|(U,V)|^{3/2}}{3/2} + \frac{|(\Phi,\Psi)|^3}{3}\bigg)\,.
\end{equation}

Combine \eqref{B1dI} and \eqref{bsverhu}. We obtain \eqref{sqrt32}.

To find the function with \eqref{pr1} and \eqref{pr2} we need the next section.

\section{Function $B=\frac{2}{9} ( y_{11}^2 + y_{12}^2) + 3 (y_{21}^2 + y_{22}^2)^{1/2})^{3/2} +\frac{2}{9} (( y_{11}^2 + y_{12}^2))^{3/2}$.}
\label{f}

It is useful if the reader thinks that $y_{11}, y_{12}, y_{21}, y_{22}$ are correspondingly $\Phi, \Psi, U, V$.

Also in what follows $dy_{11}, dy_{12}, dy_{21}, dy_{22}$ can be viewed as $\phi_1, \psi_1, u_1, v_1$ and $\phi_2, \psi_2, u_2, v_2$.

\bigskip

Let $\Bel_{n+m}(x)$ be a real-valued function of $n+m$ variables
$x=(x_1,\dots,x_n,x_{n+1},\dots,x_{n+m})$. Define a function $\Bel_{nk+m}(y)$
of $n$ vector valued variables $y_i=(y_{i1},\dots,y_{ik})$, $1\le i\le n$, and
$m$ scalar variables $y_i$, $n+1\le i\le n+m$, as follows:
$$
\Bel_{nk+m}(y)=\Bel_{n+m}(x),
$$
where
\begin{align*}
x_i&=\|y_i\|:=\Big(\sum_{j=1}^k y_{ij}^2\Big)^\half\qquad&\text{for }i\le n,
\\
x_i&=y_i\qquad&\text{for }i>n.
\end{align*}
Omitting indices we shall denote by $\frac{d^2\Bel}{dx^2}$ and
$\frac{d^2\Bel}{dy^2}$ the Hessian matrices of $\Bel_{n+m}(x)$ and
$\Bel_{nk+m}(x)$ respectively.

\section{Hessian of a vector-valued function}
\label{hessian}

\begin{lemma}
\label{Hess1}
Let $P_j$ be the following operator from $\R^k$ to $\R$:
$$
P_jh=\frac{(h,y_j)}{x_j},
$$
i.e., it gives the projection to the direction $y_j$, and let $P$ be the
block-diagonal operator from
$\R^{kn+m}=\R^k\oplus\R^k\oplus\dots\oplus\R^k\oplus\R\oplus\dots\oplus\R$ to
$\R^{n+m}=\R\oplus\R\oplus\dots\oplus\R\oplus\R\oplus\dots\oplus\R$ whose first
$n$ diagonal elements are $P_j$ and the rest is identity. Then
$$
\frac{d^2\Bel}{dy^2}=P^*\frac{d^2\Bel}{dx^2}P+\text{diag}\,\left\{(I-P_i^*P_i)\frac1{x_i}
\frac{\partial\Bel}{\partial x_i}\right\},
$$
or
\begin{align*}
d^2\Bel &= \sum_{i,j=1}^n\frac{\partial^2\Bel}{\partial x_i\partial
x_j}\cdot\frac{\sum_{s=1}^ky_{is}dy_{is}}{x_i}\cdot\frac{\sum_{r=1}^ky_{jr}dy_{jr}}{x_j}
\\
&+2\sum_{i=1}^n\sum_{j=n+1}^{n+m}\frac{\partial^2\Bel}{\partial x_i\partial
x_j}\cdot\frac{\sum_{s=1}^ky_{is}dy_{is}}{x_i}\cdot dy_j
\\
&+\sum_{i=n+1}^{n+m}\sum_{j=n+1}^{n+m}\frac{\partial^2\Bel}{\partial
x_i\partial x_j}\cdot dy_i\cdot dy_j
\\
&+\sum_{i=1}^n\frac1{x_i}\frac{\partial\Bel}{\partial x_i}
\cdot\left(\sum_{j=i}^kdy_{ij}^2-\Big(\frac{\sum_{j=1}^ky_{ij}dy_{ij}}{x_i}\Big)^2\right)\,.
\end{align*}
\end{lemma}

\subsection{Positive definite quadratic forms}
\label{qf}

Let
\begin{equation*}
Q=Ax^2+2Bxy+Cy^2
\end{equation*}
be a positive definite quadratic form. We are interested in the best possible
constant $D$ such that
\begin{equation*}
Q\ge 2D|x|\,|y|\qquad\text{for all } x,y\in\R\,.
\end{equation*}
After dividing this inequality over $|x|\,|y|$ we get
$$
At\pm2B+\frac{C}t\ge 2D\qquad\text{for all } t\in\R\setminus\{0\}\,,
$$
The left-hand side has its minimum at the point $t=\sqrt{\frac CA}$\,.
Therefore the best $D$ is $\sqrt{AC}-|B|$.

Now we would like to present $Q$ as a sum of three squares:
$$
Q=D(\tau x^2+\frac1\tau y^2)+(\alpha x+\beta y)^2\,,
$$
which immediately implies the required estimate. By the assumption
$$
(A-D\tau)x^2+2Bxy+(C-\frac D\tau)y^2
$$
is a complete square, whence
$$
(A-D\tau)(C-\frac D\tau)=B^2
$$
or
\begin{gather*}
CD\tau^2-(AC-B^2+D^2)\tau+AD=0\,,
\\
C\tau^2-2\sqrt{AC}\tau+A=0\,.
\end{gather*}
Therefore, $\tau=\sqrt{\frac AC}$ and
\begin{equation}\label{Quadratic}
Q=(\sqrt{AC}-|B|)\Big(\sqrt{\frac AC}x^2+\sqrt{\frac CA}y^2\Big)+
|B|\sqrt{\frac AC}\Big(x+\sign B\sqrt{\frac CA}y\Big)^2
\end{equation}

\subsection{Example}
\label{example}
Let
\begin{equation}
\label{Bel3}
\Bel_2(x) = \frac{2}{9}(x_1^2+3x_2)^{3/2}+\frac{2}{9} x_1^3\,,
\end{equation}
$$
\Bel_4(y) = B_2(x);\qquad x_i=\sqrt{y_{i1}^2+y_{i2}^2}\,.
$$
Calculate the derivatives:
\begin{gather*}
\frac{\partial\Bel_2}{\partial x_1} =\frac{2}{3} x_1(\sqrt{x_1^2+3x_2}+x_1)\,,
\\
\frac{\partial\Bel_2}{\partial x_2} =\sqrt{x_1^2+3x_2}\,,
\\
A=\frac{\partial^2\Bel_2}{\partial x_1^2}
=\frac{2(\sqrt{x_1^2+3x_2}+x_1)^2}{3\sqrt{x_1^2+3x_2}}\,,
\\
B=\frac{\partial^2\Bel_2}{\partial x_1\partial x_2}
=\frac{x_1}{\sqrt{x_1^2+3x_2}}\,,
\\
C=\frac{\partial^2\Bel_2}{\partial x_2^2} =\frac3{2\sqrt{x_1^2+3x_2}}\,,
\\
D=\sqrt{AC}-|B|=1\,,
\end{gather*}
Also
\begin{equation}
\label{tau0}
\tau=\sqrt{\frac AC}=\frac{2}{3}(\sqrt{x_1^2+3x_2}+x_1)\,,
 \end{equation}
\begin{equation}
\label{1tau0}
\frac1\tau=\frac{\sqrt{x_1^2+3x_2}-x_1}{2x_2}\,.
\end{equation}

After substitution in the expressions of the preceding sections we get
\begin{align*}
d^2 \Bel_4 &= \tau\Big(\frac{y_{11}dy_{11}+y_{12}dy_{12}}{x_1}\Big)^2 +
\frac1\tau \Big(\frac{y_{21}dy_{21}+y_{22}dy_{22}}{x_2}\Big)
\\
&+\frac{\tau x_1}{\sqrt{x_1^2+3x_2}}
\Big[\frac{y_{11}dy_{11}+y_{12}dy_{12}}{x_1}+
\frac1\tau\frac{y_{21}dy_{21}+y_{22}dy_{22}}{x_2}\Big]^2
\\
&+\frac{2}{3} \Big(\sqrt{x_1^2+3x_2}+x_1\Big)
\Big(\frac{y_{12}dy_{11}-y_{11}dy_{12}}{x_1}\Big)^2
\\
&+\frac{\sqrt{x_1^2+3x_2}}{x_2}\Big(\frac{y_{22}dy_{21}-y_{21}dy_{22}}{x_2}\Big)^2
\\
&=\tau(dy_{11}^2+dy_{12}^2) + \frac1\tau (dy_{21}^2+dy_{22}^2)
+\frac{3\tau}{4x_2} \Big(\frac{y_{22}dy_{21}-y_{21}dy_{22}}{x_2}\Big)^2
\\
&+\frac{\tau x_1}{\sqrt{x_1^2+3x_2}}
\Big[\frac{y_{11}dy_{11}+y_{12}dy_{12}}{x_1}+
\frac1\tau\frac{y_{21}dy_{21}+y_{22}dy_{22}}{x_2}\Big]^2.
\end{align*}

\subsection{Verifying \eqref{tau}.}
\label{tausection}

Here using \eqref{tau0}, \eqref{1tau0} we get
\begin{equation}
\label{tau1}
\tau =\frac{2}{3}((y_{11}^2 +y_{12}^2)+3(y_{21}^2 +y_{22}^2)^{1/2})^{1/2} +(y_{11}^2 +y_{12}^2)^{1/2})\,.
\end{equation}
And henceforth
\begin{equation}
\label{1tau}
\frac1{\tau} = \frac{((y_{11}^2 +y_{12}^2)+3(y_{21}^2 +y_{22}^2)^{1/2})^{1/2} - (y_{11}^2 +y_{12}^2)^{1/2}}{2(y_{21}^2 +y_{22}^2)^{1/2}}\,.
\end{equation}

Let us now (when we know $\tau$) check the condition \eqref{tau}:
$$
\frac{3}{4}\frac{\tau}{(y_{21}^2 +y_{22}^2)^{1/2}} +\frac2{\tau} = 
\frac1{2} \frac{(y_{11}^2 +y_{12}^2)+3(y_{21}^2 +y_{22}^2)^{1/2})^{1/2} + (y_{11}^2 +y_{12}^2)^{1/2}}{(y_{21}^2 +y_{22}^2)^{1/2}} +
$$
$$
\frac{(y_{11}^2 +y_{12}^2)+3(y_{21}^2 +y_{22}^2)^{1/2})^{1/2} -(y_{11}^2 +y_{12}^2)^{1/2}}{(y_{21}^2 +y_{22}^2)^{1/2}}=
$$
$$
\frac{3(y_{11}^2 +y_{12}^2)+3(y_{21}^2 +y_{22}^2)^{1/2})^{1/2} -(y_{11}^2 +y_{12}^2)^{1/2}}{2(y_{21}^2 +y_{22}^2)^{1/2}} \ge \frac{3}{\tau}\,.
$$


\section{Explanation. Pogorelov's theorem.}
\label{expl}

We owe the reader  the explanation, where we got this function $B$, which played such a prominent part above.

We want to find a function satisfying the following propperties (in what follows $p\ge 2$):

\begin{itemize}
\item 1) $B$ is defined in the whole plane $\R^2$ and $B(u,v) = B(-u,v)= B(u,-v)$;
\item 2) $0\le B(u,v) \le (p-1) (\frac1p |u|^p + \frac1q |v|^q)$;
\item 3) Everywhere we have inequality for Hessian quadratic form $d^2 B(u,v) \ge 2|du||dv|$;
\item 4) Homogenuity: $B(c^{1/p}u,c^{1/q}v)= c\,B(u,v)$, $c>0$;
\item 5) Function $B$ should be the ``best" one satisfying 1), 2), 3).
\end{itemize}

The last statement we will understand in the following sense: $B$ must saturate inequalities to make them equalities on a natural subset of $\R^2$ in 2) and on a natural subset of the tangent bundle of $\R^2$ in 3).

Let us start with 3). Inequalities just mean that $d^2 B(u,v) \ge 2dudv\,,\, d^2 B(u,v) \ge -2dudv$ for any
$(u,v)\in \R^2$ and for any $(du,dv)\in \R^2$. In other words this is just positive definitness of matrices
\begin{equation}
\label{posdef}
\begin{bmatrix} B_{uu}, & B_{uv} -1\\
B_{vu} -1, & B_{vv}\end{bmatrix}\ge 0\,,\,\,\begin{bmatrix} B_{uu}, & B_{uv} +1\\
B_{vu} +1, & B_{vv}\end{bmatrix}\ge 0\,.
\end{equation}

Now we want that \eqref{posdef} barely occurs. In other words we want that for any $(u,v)$ one of the matrix in \eqref{posdef} would have a zero determinant.

Notice that symmetry 1) allows us to consider $B$ only in the first quadrant. Here we will assume the first matrix in \eqref{posdef} to have zero determinant in the first quadrant. 

So let us assume for $u>0, v>0$
\begin{equation}
\label{ma1}
 5)\,\,\,\,\,\,\,\,\,\,\,\,\,\,\,\,\,\,\,\,\,\,\det\begin{bmatrix} B_{uu}, & B_{uv} +1\\
B_{vu} -1, & B_{vv}\end{bmatrix}= 0\,.
\end{equation}

Let us introduce
$$
A(u,v) :=B(u,v)+uv\,.
$$
So we require
\begin{equation}
\label{ma2}
 \det\begin{bmatrix} A_{uu}, & A_{uv} \\
A_{vu} , & A_{vv}\end{bmatrix}= 0\,.
\end{equation}

Returning to saturation of 2): we require that $B(u,v) = \phi(u,v):=(p-1)(\frac1p u^p + \frac1q v^q)$ at a non-zero point. By homogenuity 4)  we have this equality on a whole curve $\Gamma$ invariant under transformations $u\rightarrow c^{1/p} u, v\rightarrow c^{1/q}v$.

\begin{equation}
\label{gamma}
B(u,v) = \phi(u,v):=(p-1)(\frac1p u^p + \frac1q v^q)\,\,\,\text{on the curve}\,\,\, v^q = \gamma u^p\,.
\end{equation}

Notice that $\gamma$ is unknown at this moment. We are going to solve \eqref{ma2}, \eqref{gamma} in the sense that our solution satisfies \eqref{posdef}, 1), 2), 3), 4).

\medskip

\noindent{\bf Remark}. We strongly suspect that the solution like that is still non-unique.
On the other hand one cannot ``improve" 1), 2), 3), 4) by, say, changing $2$ in 3) to a bigger constant, or making a constant $p-1$ in 2) smaller.

\medskip

Recall that we have also the symmetry conditions on $A(u,v) +uv=: B(u,v)$. They are 
$$
B(-u,v)= B(u,v)\,,\,\, B(u,-v) = B(u, v)\,.
$$
We assume the smoothness of $B$. It is a little bit ad hoc assumption, and we will be using it as such, namely, we will assume it when it is convenient and we will be on guard not to come to a contradiction. Anyway, assuming now the smoothness of $B$ on the $v$-axis we get that the symmetry implies the Neumann boundary condition on $B$ on $v$-axis: $ \frac{\partial}{\partial u} B(0,v) =0$, that is
\begin{equation}
\label{neu}
\frac{\partial}{\partial u} A(0,v) = v\,.
\end{equation}

Solving the homogeneous Monge-Amp\`ere equation is the same as building a surface of zero gaussian curvature. We base the following on a Theorem of Pogorelov \cite{Pog}. The reader can see the algorithm in \cite{VaVo2}. So we will be brief. Solution $A$ must have the form
\begin{equation}
\label{p1}
A(u,v) = t_1 \cdot u + t_2\cdot v -t\,,
\end{equation}
where $t_1:= A_u(u,v), t_2:=A_v(u,v), t(u,v)$ are unknown function of $u,v$, but, say, $t_1, t_2$ are certain functions of $t$. Moreover, Pogorelov's theorem says that
\begin{equation}
\label{p2}
u\cdot dt_1 + v\cdot dt_2 -dt=0\,,\,\,\text{meaning}\,\,u\cdot \frac{dt_1}{dt} + v\cdot \frac{dt_2}{dt} -1=0\,.
\end{equation}

We write homogenuity condition 4) as follows $A(c^{1/p}u, c^{1/q}v) =cA(u,v)$, differentiate in $c$ and plug $c=1$. Then we obtain
\begin{equation}
\label{p3}
A(u,v) = \frac{1}{p}t_1 \cdot u + \frac{1}{q}t_2\cdot v \,,
\end{equation}
which being combined with \eqref{p1} gives
\begin{equation}
\label{p4}
\frac1q t_1 \cdot u + \frac1p t_2\cdot v -t=0\,.
\end{equation}
Notice a simple thing, when $t$ is fixed \eqref{p2} gives us the equation of a line in $(u,v)$ plane. Call this line $L_t$. Functions $t_1, t_2$ are certain (unknown at this moment) functions of $t$, so again, for a fixed $t$ equation \eqref{p4} also gives us a line. Of course this must be $L_t$. Comparing the coefficients we obtain differential equations on $t_1, t_2$:

\begin{equation}
\label{p5}
q\frac{d t_1}{t_1} = \frac{dt}{t}\,,\,p\frac{d t_2}{t_2} = \frac{dt}{t}\,.
\end{equation}
We write immediately the solutions in the following form:
\begin{equation}
\label{p6}
t_1(t)= pC_1|t|^{\frac1q}\,,\,t_2(t)= qC_2 |t|^{\frac1p}\,.
\end{equation}
Plugging this into \eqref{p3} one gets
\begin{equation}
\label{A=}
A(u,v) = C_1t^{\frac1q} u + C_2 t^{\frac1p} v \,,\,\,B(u,v) = C_1t^{\frac1q} u + C_2 t^{\frac1p} v + uv\,,
\end{equation}
where  $t(u,v)$ (see \eqref{p4}) is defined from the following implicit formula

\begin{equation}
\label{t=}
t=\frac{p}q C_1 t^{\frac1q} u + \frac{q}p C_2 t^{\frac1p} v\,.
\end{equation}

To define unknown constants $C_1, C_2$ we have only one boundary condition \eqref{neu}.  However we have one more condition. It is a free boundary condition (we think that $p\ge 2\ge q$)
\begin{equation}
\label{fbc1}
B(u,v) = \phi(u,v):= (p-1)(\frac1p u^p + \frac1q v^q)\,\,\text{on the curve}\,\,\Gamma:=\{v^q= \gamma^q u^p\}\,.
\end{equation}
This seems to be not saving us because we have three unknowns $C_1, C_2, \gamma$ and two conditions: \eqref{neu} and {\eqref{fbc1}. But we will require in addition that $B(u,v)$ and $\phi(u,v)$ have the same tangent plane on the curve $\Gamma$:
\begin{equation}
\label{fbc2}
\frac{B_u(u,v)}{B_v(u,v)} = \frac{\phi_u(u,v)}{\phi_v(u,v)}\,\,\text{on the curve}\,\,\Gamma =\{v^q= \gamma^q u^p\}\,.
\end{equation}

Now we are going to solve \eqref{neu}, \eqref{fbc1}, \eqref{fbc2}, to find $C_1, C_2, \gamma$ and plug them into \eqref{A=} and \eqref{t=}. 

First of all
$$
v= A_u(0,v) = t_1(0,v)\,.
$$
So $ v/pC_1 = t(0,v)^{\frac1q}$ from \eqref{p6}. Plug $u=0$ into \eqref{t=} to get $t(0,v)^{\frac1q} =\frac{q}p C_2 v$. Combining we get
$$
C_1C_2 =\frac1q\,.
$$

Now we use \eqref{fbc2}. 
$$
\frac{t_1 -v}{t_2-u} = \frac{u^{p-1}}{v^{q-1}}= \frac{u^p}{v^q}\frac{v}{u}= \frac1{\gamma^q}\frac{v}{u}\,.
$$
Using \eqref{p6} we get
\begin{equation}
\label{G2}
\frac{pC_1 t^{\frac1q} -v}{qC_2t^{\frac1p}-u} = \frac1{\gamma^q}\frac{v}{u}\,.
\end{equation}
Let us write $\Gamma$ as $ u^p =\frac1{\gamma} uv$ or $v^q= \gamma^{q-1} uv$, and let us write on $\Gamma$
\begin{equation}
\label{ab}\begin{cases}
t^{\frac1q} = av\\
t^{\frac1p} = bu
\end{cases}\end{equation}
The reader will easily see from what follows that $a, b$ are constants.  From \eqref{G2}
\begin{equation}
\label{G21}
(pC_1a -1)\gamma^q uv =(qC_2b-1)uv \,.
\end{equation}
Also from \eqref{ab}
\begin{equation}
\label{G3}
\frac{a^q}{b^p}=\frac1{\gamma^q}\,,
\end{equation}
and from \eqref{ab}  and \eqref{t=}
\begin{equation}
\label{G4}
ab = \frac{p}q C_1 a+ \frac{q}p C_2 b\,.
\end{equation}
From \eqref{ab}, \eqref{fbc1} it follows
\begin{equation}
\label{G5}
C_1 a+  C_2 b- 1=(p-1) (\frac1p\cdot \frac1{\gamma} + \frac1q\cdot \gamma^{q-1})\,.
\end{equation}
We already proved 
\begin{equation}
\label{G6}
C_1 C_2 =\frac1q\,.
\end{equation}
We have five equations \eqref{G2}--\eqref{G6} on five unknowns $C_1, C_2, a, b, \gamma$.

One solution is obvious: 
$$
\gamma=1\,, a=qC_2\,, b=pC_1\,, p^p C_1^p = q^q C_2^q\,,
$$
from where one finds
\begin{equation}
\label{CC}
C_1 = \frac1p p^{\frac1p}\,,\, C_2 = \frac1q p^{\frac1q}\,.
\end{equation}
Therefore,
\begin{equation}
\label{Buv}
B(u,v) = \frac1p p^{\frac1p} t^{\frac1q} u + \frac1q p^{\frac1q} t^{\frac1p} v -uv\,,
\end{equation}
where $t$ is defined from
\begin{equation}
\label{tuv}
t = \frac1q p^{\frac1p} t^{\frac1q} u + \frac1p p^{\frac1q} t^{\frac1p} v\,.
\end{equation}

\medskip 

If we specify $p=3, q=\frac32$ we get
\begin{equation}
\label{CC3}
C_1=\frac13 3^{\frac13}\,,\, C_2 = \frac23 3^{\frac23}\,.
\end{equation}
\begin{equation}
\label{tuv3}
t^{\frac23} = \frac23 3^{\frac13} t^{\frac13} u + \frac13 3^{\frac23}  v\,,
\end{equation}
and solving the quadratic equation on $s:=t^{\frac13}$: $ s^2 -2C_1 us -\frac{C_2}{2} v =0$, we get (the right root will be with $+$ sign)
\begin{equation}
\label{s}
t^{\frac13}(u,v)=s= C_1 u + \sqrt{C_1^2 u^2 + \frac{C_2}{2} v}\,.
\end{equation}
Therefore, $B(u,v)$ being equal to $C_1 s^2 u + C_2 s v -uv$  is ($C_1C_2= \frac23$, see \eqref{G6})
$$
B(u,v) = C_1 u (2C_1us + \frac{C_2}2 v) + C_2 vs -uv= (2C_1^2 u^2  +C_2v) s +\frac12C_1C_2 uv-uv\,,
$$
and so
$$
B(u,v) =(2C_1^2 u^2  +C_2v) ( C_1 u +\sqrt{C_1^2 u^2 + \frac{C_2}{2} v}) -\frac23 uv\,,
$$
$$
=(2C_1^2 u^2  +C_2v)\sqrt{C_1^2 u^2 + \frac{C_2}{2} v} + 2C_1^3 u^3 + (C_1C_2-\frac23) uv\,.
$$
The last term disappears (see \eqref{G6}), and we get
$$
B(u,v) = 2(C_1^2 u^2 + \frac{C_2}{2} v)\sqrt{C_1^2 u^2 + \frac{C_2}{2} v} + 2C_1 u^3=
2C_1^3 ( u^2 + \frac{C_2}{2C_1} v)^{\frac32} + 2C_1^3 u^3\,.
$$
Finally from \eqref{CC3}
\begin{equation}
\label{Buv3}
B(u,v)= \frac29 ((u^2 + 3v)^{\frac32} +u^3)\,.
\end{equation}

This is exactly the function in \eqref{Bel3}. This function gave us our main theorem for $p=3$. We have just explained how we got it.

\medskip

By the way, in this particular case the transcendental equation on $\gamma$ becomes the usual cubic equation on $\sqrt{\gamma}$: $ 2\sqrt{\gamma} +1 =4-\frac1{\gamma}$, which has only one real solutions $\gamma=1$.


\bigskip

\section{Explanation. Pogorelov's theorem again.}
\label{expl}

We owe the reader  the explanation, why we chose the function $A(u,v)=B(u,v)+uv$ rather than $A(u,v)=B(u,v)-uv$ to have the degenerate Hessian form.

We want to find a function satisfying the following properties (in what follows $p\ge 2$):

\begin{itemize}
\item 1) $B$ is defined in the whole plane $\R^2$ and $B(u,v) = B(-u,v)= B(u,-v)$;
\item 2) $0\le B(u,v) \le \phi(u,v)=(p-1) (\frac1p |u|^p + \frac1q |v|^q)$;
\item 3) Everywhere we have inequality for Hessian quadratic form $d^2 B(u,v) \ge 2|du||dv|$;
\item 4) Homogenuity: $B(c^{1/p}u,c^{1/q}v)= c\,B(u,v)$, $c>0$;
\item 5) Function $B$ should be the ``best" one satisfying 1), 2), 3).
\end{itemize}

\begin{enumerate}
\item
What do we mean by $best$ function? We would like $B$ to be the 'largest' function below $\phi(u,v)$ such that the convexity condition in 3) holds. We expect that such a function should equal the upper bound $\phi(u,v)$ at some point(s) and the inequality in 3) should be equality where possible. 
\item
Due to the symmetry in 1), we can restrict our attention to $\{u>0, v>0\}$.  
\item If we have at some $(u,v)$, $B(u,v)= \phi(u,v)$, then condition 4) implies that $B(c^{1/p}u, c^{1/q}v)= cB(u,v)=c\phi(u,v)=\phi(c^{1/p}u, c^{1/q}v)$. Hence they remain equal on a curve $\{(u,v): v^q = \gamma^q u^p\}$ for some $\gamma$.
\item The condition $<d^2B\cdot(u,v),(u,v)> \geq 2|u||v|$ means that the 'directional convexity' in direction $(u,v)$ stays above the value $2|u||v|$. This means that the directional convexity of $B$ is above that of both the functions $uv$ and $-uv$. Equivalently we are asserting the positive definiteness of the matrices:

\begin{equation}\label{N1}
\begin{pmatrix} B_{uu} & B_{uv}-1 \\ B_{vu}-1 & B_{vv}\end{pmatrix}\geq 0, \begin{pmatrix} B_{uu} & B_{uv}+1 \\ B_{vu}+1 & B_{vv}\end{pmatrix}\geq 0.
\end{equation}

\item In order to optimize (\eqref{N1}), we require that one of the matrices is degenerate (with $"=0"'$). Suppose that the first matrix is degenerate. This means that the function $A(u,v) = B(u,v)-uv$ has a degenerate Hessian. At every point, one of its two non-negative eigenvalues is $0$, and the function has $0$ convexity in the corresponding eigendirection. Since the matrix is positive definite, it follows that $0$ is the minimal eigenvalue, hence the graph of this function is a surface with gaussian curvature $0$.

Moreover the directional convexity of $B - uv$ is greater than that of $B+uv$ in directions of negative slope and less than in directions of positive slope. If we want $B+uv$ to have non-degenerate positive Hessian, then the degeneracy of $B-uv$ must occur in the positive slope direction. 
\end{enumerate}

Let us analyse the function $A(u,v)= B(u,v)-uv$. A theorem of Pogorolev tells us that $A$ will be a linear function on lines of degeneracy. That is, it will have the form:
\begin{equation}\label{N2}
A(u,v) = t_1u+t_2v-t
\end{equation}
where $t_1(u,v)$, $t_2(u,v)$ and $t(u,v)$ are constant on the lines given by 
\begin{equation}\label{N3}
\frac{dt_1}{dt}u+\frac{dt_2}{dt}v-1 = 0.
\end{equation}
We can say two things about the coefficient functions, that the eigenlines that intersect the positive $y$ axis must also have $\frac{dt_1}{dt_2}\leq 0$ and $\frac{dt_2}{dt}\geq 0$ - this information comes from (\eqref{N3}) and the fact that the eigenlines have positive slope. At the moment we know nothing else about the coefficient functions. We will use the various boundary conditions on $B$, hence on $A$ to determine them. 

\begin{enumerate}
\item First observe that since $B(u,v)= B(-u,v)=B(u,-v)$, we may expect that $B$ is smooth on at least one of the two axes, assume on the $y$ axis, and hence the corresponding derivative $\partial_u B(0,v) = 0$. This means:
\begin{equation}\label{N4}
\partial_u A(0,v) = -v.
\end{equation}
\item We already assumed that 
\begin{equation}\label{N5}
B(u,v) = \phi(u,v) = (p-1)(\frac{u^p}{p} + \frac{v^q}{q})
\end{equation} on some curve $\Gamma=\{v^q=\gamma^qu^p\}$.
\item Let us also assume that the tangent planes of $B$ and $\phi$ agree on $\Gamma$. This means that the gradients of the two functions $B(u,v)-z$ and $\phi(u,v)-z$ should be parallel at the points $(u,v,\phi(u,v))$ where $(u,v)\in \Gamma$. Therefore
\[(\partial_u\phi, \partial_v \phi, -1) = \lambda (\partial_u B, \partial_v B, -1),\] which implies $\lambda = 1$ and 
\begin{equation}\label{N6}
B_u(u,v) = (p-1)u^{p-1}, B_v(u,v) = (p-1)v^{q-1}
\end{equation}
on the curve $\Gamma$. Similarly on $\Gamma$,
\begin{equation}\label{N7}
A_u(u,v) = (p-1)u^{p-1}-v, A_v(u,v) = (p-1)v^{q-1}-u.
\end{equation}
\end{enumerate} 

Recall:
\begin{equation}\label{N8}
A(u,v) = t_1u+t_2v-t
\end{equation}
where $t_1(u,v)=A_u(u,v)$, $t_2(u,v)=A_v(u,v)$ and $t(u,v)$ are constant on the lines given by 
\begin{equation}\label{N9}
\frac{dt_1}{dt}u+\frac{dt_2}{dt}v-1 = 0.
\end{equation}
We also have the homogeneity condition: $A(c^{1/p}u,c^{1/q}v)=cA(u,v)$. Differentiating this with respect to $c$ and setting $c=1$ gives:
\begin{eqnarray}
A(u,v)&=& \frac{1}{p}A_u(u,v)u+\frac{1}{q}A_v(u,v)v \label{N10}\\
&=& \frac{1}{p}t_1 u+\frac{1}{q}t_2 v.\label{N11}
\end{eqnarray}
Comparing (\eqref{N8}) and (\eqref{N11}), we have 
\begin{equation}\label{N12}
\frac{1}{q}t_1u+\frac{1}{p}t_2v -t = 0.
\end{equation}
Now comparing (\eqref{N9}) and (\eqref{N12}) gives
\begin{equation}\label{N13}
\frac{dt_1}{dt}=\frac{1}{q}\frac{t_1}{t}, \frac{dt_2}{dt}=\frac{1}{p}\frac{t_2}{t}.
\end{equation}
Solving these differential equations, we have 
\begin{equation}\label{N14}
t_1(t)= C_1 |t|^{1/q}, t_2(t)= C_2 |t|^{1/p}.
\end{equation}

Putting this into (\eqref{N12}) gives:
\begin{equation}\label{N15}
t = \frac{1}{q}C_1|t|^{1/q} u + \frac{1}{p}C_2 |t|^{1/p}v
\end{equation}

Let us make two observations:
Recall that if our eigenline intercepts the positive $y$ axis and has positive slope, then $\frac{dt_1}{dt_2}=\frac{q}{p}\frac{C_1}{C_2}|t|^{\frac{1}{q}-\frac{1}{p}}\leq 0$ and $\frac{dt_2}{dt}\geq 0$. If $t>0$, then
$\frac{dt_2}{dt} = \frac{1}{p}C_2|t|^{-1/q}$, and if $t<0$, then 
$\frac{dt_2}{dt} = -\frac{1}{p}C_2|t|^{-1/q}$. We conclude from this:

\begin{enumerate}
\item If $t>0$, then $C_1C_2 \leq 0$ and $C_2\geq 0$, hence $C_1\leq 0$,
\item If $t<0$, then $C_1C_2 \leq 0$ and $C_2\leq 0$, hence $C_1\geq 0$.
\end{enumerate}

Let us bring in the following: $t_1= A_u(0,v)=-v$. The first equality is from Pogorolev and the second is the boundary condition \eqref{N4}. Then \eqref{N14} implies that 
\begin{equation}\label{N16}
-v = C_1 |t(0,v)|^{1/q}
\end{equation} and \eqref{N15} implies that 
\begin{equation}\label{N17}
t(0,v) = \frac{1}{p}C_2 |t(0,v)|^{1/p}v.
\end{equation}

Conclude:
\begin{enumerate}
\item If $v>0$, then $C_1<0$. The previous observations imply $t>0$ and $C_2\geq 0$. We are concerned at present with this case of positive $y$ intercept.
\item From \eqref{N16} and \eqref{N17}, we conclude
\begin{equation}\label{N18}
C_1C_2 = - p.
\end{equation}
\end{enumerate}

Next from \eqref{N7}, we know that on $\Gamma$,
\begin{equation}\label{N19}
t_1 = (p-1)u^{p-1}-v = (\frac{p-1}{\gamma}-1)v, t_2 = (p-1)v^{q-1}-u = ((p-1)\gamma^{q-1} - 1)u.
\end{equation}
In terms of $t$, this says
\begin{equation}\label{N20}
C_1 t^{1/q} = (\frac{p-1}{\gamma}-1)v, C_2 t^{1/p} = ((p-1)\gamma^{q-1} - 1)u.
\end{equation}
Write on $\Gamma$
\begin{equation}
\label{N21}\begin{cases}
t^{\frac1q} = aC_2v\\
t^{\frac1p} = bC_1u
\end{cases}\end{equation}
Note that $a\geq 0$ and $b\leq 0$ due to the signs of $C_1$ and $C_2$.
Substituting in \eqref{N20} and using \eqref{N18} gives
\begin{equation}\label{N22}
 a = \frac{1}{p}-\frac{1}{q\gamma}, b = \frac{1}{p}-\frac{1}{q}\gamma^{q-1}.
\end{equation}
Note that \eqref{N21} also implies that 
\begin{equation}\label{N23}
\frac{a^q C_2^q}{|b|^p|C_1|^p} = \frac{1}{\gamma^q}
\end{equation}
Hence \eqref{N22} and \eqref{N23} imply
\begin{equation}\label{N24}
(\frac{\gamma}{p}-\frac{1}{q})^q C_2^q = (\frac{1}{q}\gamma^{q-1}-\frac{1}{p})^p |C_1|^p.
\end{equation}
\eqref{N24}, \eqref{N18} and the fact $pq=p+q$ imply that
\begin{equation}\label{N25}
C_2 = \left(\frac{p((p-1)\gamma^{\frac{1}{p-1}}-1)^{p-1}}{\gamma-(p-1)}\right)^{\frac{1}{p}}.
\end{equation}

Next observe that \eqref{N15}, \eqref{N18} and \eqref{N21} imply 
\begin{equation}\label{N26}
ab = \frac{1}{q}a+\frac{1}{p}b
\end{equation}
and hence by \eqref{N22}
\begin{equation}\label{N27}
(\frac{1}{p}-\frac{1}{q\gamma})(\frac{1}{p}-\frac{1}{q}\gamma^{q-1}) = 
\frac{1}{pq}+\frac{1}{p^2}-\frac{1}{q^2\gamma}-\frac{1}{pq}\gamma^{q-1}
\end{equation}
The equation that follows from making substitutions into the boundary condition \eqref{N5} $B=\phi$ on $\Gamma$ and $A=B-uv$ gives no new relationship. So we can avoid its consideration.

Simplifying \eqref{N27} shows that $\gamma$ is solution to the equation
\begin{equation}\label{N28}
\gamma^{q-1}-(q-1)\gamma+ 2-q = 0.
\end{equation}

The rest of the analysis is yet to be done. However note that $\frac{B_u}{u}=\frac{\phi_u}{u}$ on $\Gamma$, and on the 
corresponding eigenline, we can understand it by using the fact that $A_u = B_u -v$ is constant. This may help later.

\section{The case when $p=3$ and $q=\frac{3}{2}$}
Observe that by setting $\delta=\gamma^{q-1}$, we can rewrite \eqref{N28} as
\begin{equation}\label{P1}
\delta^{p-1}-(p-1)\delta+2-p = 0.
\end{equation}
Let us analyse the case when $p=3$. Then this equation becomes
\begin{equation}\label{P2}
\delta^2-2\delta-1=0
\end{equation}
whose unique positive solution is $\delta= 1+\sqrt{2}$. Therefore 
\begin{equation}\label{P2}
\gamma = (1+\sqrt{2})^2= 3+2\sqrt{2}.
\end{equation}
Then using \eqref{N18}, \eqref{N22} and \eqref{N25}, we obtain
\begin{equation}\label{P3}
a= -\frac{5}{3} + \frac{4\sqrt{2}}{3}, b = \frac{1}{3}-\frac{2}{3}(3+2\sqrt{2})^{1/2}
\end{equation} and
\begin{equation}\label{P4}
C_1 = \frac{-3^{\frac{2}{3}}(1+2\sqrt{2})^{1/3}}{(2\sqrt{3+2\sqrt{2}}-1)^{2/3}}, C_2 = \frac{3^{1/3}(2\sqrt{3+2\sqrt{2}}-1)^{2/3}}{(1+2\sqrt{2})^{1/3}}
\end{equation}
Now we will explicitly find $B(u,v)$. Recall
\begin{eqnarray*}
B(u,v)&=& \frac{1}{p}t_1u+\frac{1}{q}t_2v+uv \\
&=& \frac{C_1}{p}t^{1/q}u+\frac{C_2}{q}t^{1/p}v+uv\\
&=& \frac{C_1}{3}t^{2/3}u+\frac{2C_2}{3}t^{1/3}v+uv.
\end{eqnarray*}
\begin{eqnarray*}
t&=&\frac{1}{q}C_1t^{1/q}u+\frac{1}{p}C_2t^{1/p}v\\
&=&\frac{2}{3}C_1t^{2/3}u+\frac{1}{3}C_2t^{1/3}v.
\end{eqnarray*}
Let $s=t^{1/3}$. Then we have $s^2-\frac{2}{3}C_1u s-\frac{1}{3}C_2v=0$
and \[ s=\frac{C_1}{3}u+\frac{1}{3}\sqrt{C_1^2u^2+3C_2v}.\]
\begin{eqnarray*}
B(u,v)&=& \frac{C_1}{3}s^2u+\frac{2}{3}C_2sv+uv\\
&=& \frac{C_1}{3}\left(\frac{2C_1^2}{9}u^2+\frac{3C_2v}{9}+\frac{2C_1u}{9}\sqrt{C_1^2u^2+3C_2v}\right)u \\
&& + \frac{2}{9}C_1C_2uv+\frac{2}{9}C_1C_2uv+\frac{2}{9}C_2v\sqrt{C_1^2u^2+3C_2v}+uv.
\end{eqnarray*}
Use the fact that $C_1C_2=-3$ to simplify and obtain:
\begin{equation}\label{P5}
B(u,v)=\frac{2}{27}(C_1^2u^2+3C_2v)^{3/2}+\frac{2}{27}C_1^3u^3.
\end{equation}

\begin{eqnarray}
B_u &=& \frac{2}{9}C_1^3u^2+\frac{1}{9}(C_1^2u^2+3C_2v)^{1/2}2C_1^2u  \nonumber\\
B_v &=&\frac{C_2}{3}(C_1^2u^2+3C_2v)^{1/2}  \nonumber\\
B_{uu}&=&\frac{2}{9}C_1^2\left[\frac{(\sqrt{C_1^2u^2+3C_2v}+C_1u)^2}{\sqrt{C_1^2u^2+3C_2v}}\right] \nonumber\\
B_{uv}&=&\frac{|C_1|uv}{\sqrt{C_1^2u^2+3C_2v}} \nonumber\\
B_{vv}&=&\frac{C_2^2}{2\sqrt{C_1^2u^2+3C_2v}} \nonumber\\
\tau &:=&\sqrt{\frac{B_{uu}}{B_{vv}}}= \frac{2}{3}\frac{|C_1|}{C_2}(\sqrt{C_1^2u^2+3C_2v}+C_1u)\label{P6} \\
\frac{1}{\tau}&=& \frac{\sqrt{C_1^2u^2+3C_2v}-C_1u}{2|C_1|v}\label{P7} \\
\frac{B_u}{u} &=& \frac{2}{9}C_1^3u+\frac{1}{9}(C_1^2u^2+3C_2v)^{1/2}2C_1^2 \nonumber\\
\frac{B_v}{v}&=& \frac{C_2(C_1^2u^2+3C_2v)^{1/2}}{3v}.\nonumber
\end{eqnarray}
We can use  $|C_1|C_2=3$ to deduce $\frac{B_u}{u} = \tau$. 
Next we compute the quadratic form associated with $B$ by using the formulation before:
\begin{eqnarray*}
&& Q(dx,dy) = B_{uu}dx^2+2B_{uv}dxdy+B_{vv}dy^2 \\
&=& \left(\sqrt{B_{uu}B_{vv}}-|B_{uv}|\right)\left(\sqrt{\frac{B_{uu}}{B_{vv}}}dx^2+\sqrt{\frac{B_{vv}}{B_{uu}}}dy^2\right) + |B_{uv}|\sqrt{\frac{B_{uu}}{B_{vv}}}\left(dx +\sign{(B_{uv})}\sqrt{\frac{B_{vv}}{B_{uu}}}dy\right)^2 \\
&=& \left(\sqrt{B_{uu}B_{vv}}-|B_{uv}|\right)\left(\tau dx^2+\frac{1}{\tau}dy^2\right)+|B_{uv}|\tau \left(dx+\sign{(B_{uv})}\frac{1}{\tau}dy\right)^2.
\end{eqnarray*}

Now let $\Bel (y_{11},y_{12},y_{21},y_{22}):=B(\sqrt{y_{11}^2+y_{12}^2},
\sqrt{y_{11}^2+y_{12}^2}) = B(x_1, x_2)$. Then the associated quadratic
form becomes

\begin{align*}
d^2 \Bel &= \tau\Big(\frac{y_{11}dy_{11}+y_{12}dy_{12}}{x_1}\Big)^2 +
\frac1\tau \Big(\frac{y_{21}dy_{21}+y_{22}dy_{22}}{x_2}\Big)
\\
&+\frac{\tau|C_1|x_1}{\sqrt{C_1^2x_1^2+3C_2x_2}}
\Big[\frac{y_{11}dy_{11}+y_{12}dy_{12}}{x_1}+
\frac1\tau\frac{y_{21}dy_{21}+y_{22}dy_{22}}{x_2}\Big]^2
\\
&+\Big(\frac{B_u}{u}=\tau\Big)
\Big(\frac{y_{12}dy_{11}-y_{11}dy_{12}}{x_1}\Big)^2
\\
&+\Big(\frac{B_v}{v} = \frac{C_2(C_1^2x_1^2+3C_2x_2)^{1/2}}{3x_2}\Big)
\Big(\frac{y_{22}dy_{21}-y_{21}dy_{22}}{x_2}\Big)^2
\\
&=\tau(dy_{11}^2+dy_{12}^2) + \frac1\tau (dy_{21}^2+dy_{22}^2)\\
&+\Big(\frac{C_2(C_1^2x_1^2+3C_2x_2)^{1/2}}{3x_2}-\frac{1}{\tau}\Big) \Big(\frac{y_{22}dy_{21}-y_{21}dy_{22}}{x_2}\Big)^2
\\
&+\frac{\tau |C_1|x_1}{\sqrt{C_1^2x_1^2+3C_2x_2}}
\Big[\frac{y_{11}dy_{11}+y_{12}dy_{12}}{x_1}+
\frac1\tau\frac{y_{21}dy_{21}+y_{22}dy_{22}}{x_2}\Big]^2\\
&=\tau(dy_{11}^2+dy_{12}^2) + \frac1\tau (dy_{21}^2+dy_{22}^2)\\
&+\Big(\frac{3C_2\tau}{4C_1^2x_2}\Big) \Big(\frac{y_{22}dy_{21}-y_{21}dy_{22}}{x_2}\Big)^2
\\
&+\frac{\tau |C_1|x_1}{\sqrt{C_1^2x_1^2+3C_2x_2}}
\Big[\frac{y_{11}dy_{11}+y_{12}dy_{12}}{x_1}+
\frac1\tau\frac{y_{21}dy_{21}+y_{22}dy_{22}}{x_2}\Big]^2.
\end{align*}

In order for the quadratic form to have the self-improving property,
we need 
\begin{equation}\label{P8}
\frac{3C_2\tau}{4C_1^2x_2}+\frac{2}{\tau} \geq \frac{c}{\tau}
\end{equation}
for suitable constant $c$.
In fact if $\frac{C_2}{C_1^2}=1$, we know that $c=3$. This suggests that the right constant is $2+\frac{C_2}{C_1^2}\approx 3.276142375$.( 
Calculation gives $|C_1|\approx 1.329660319$ and $C_2\approx 2.256215334$, hence $\frac{C_2}{C_1^2}\approx 1.276142375$.)

If the rest of the process is the same as with the previous estimate, 
then the over all constant estimate would be approximately
\[\frac{2\sqrt{2}}{\sqrt{3.276142375}}\approx 1.562656814.\]


\bigskip

\section{The proof of Theorem \ref{less2} for general $q\in(1,2]$. The sharpness.}
\label{general}

Recall that we found for $1<q\le 2\le p<\infty\,,\,1/p+1/q=1$ the following function

\begin{equation}
\label{Bq}
B(u,v)=B_{q}(u,v)= \frac{p^{\frac1p}}{p} t^{\frac1q} u + \frac{p^{\frac1q}}{q}t^{\frac1p}v -uv\,,\,\,\text{where}
\end{equation}
\begin{equation}
\label{timpl}
t=t(u,v)\,\,\text{is the solution of}\,\,t=\frac{p^{\frac1p}}{q} t^{\frac1q} u + \frac{p^{\frac1q}}{p}t^{\frac1p}v \,.
\end{equation}
Our goal is to represent the Hessian form of this implicitely given $B$ as a sum of squares. This requires some calculations.

\begin{equation}
\label{Bu}
B_u= \frac{p^{\frac1p}}{p} t^{\frac1q}  -v + \frac1{pq} S \frac{t'_u}{t}\,,
\end{equation}
\begin{equation}
\label{Bv}
B_v= \frac{p^{\frac1q}}{q} t^{\frac1p}  -u + \frac1{pq} S \frac{t'_v}{t}\,,
\end{equation}
where
\begin{equation}
\label{skobkaS}
S:= p^{\frac1p} t^{\frac1q}u + p^{\frac{1}{q}} t^{\frac{1}{p}}v\,.
\end{equation}

Also
$$
t'_u= \frac{p^{\frac1p}}{p} t^{\frac1q}\cdot \frac{t}{t-\frac{p^{\frac1p}}{q^2} t^{\frac1q}u-\frac{p^{\frac1q}}{p^2} t^{\frac1p}v}\,,
$$
which, after using \eqref{timpl}, \eqref{skobkaS} gives
\begin{equation}
\label{tut}
\frac{t'_u}{t} = p\cdot p^{\frac1p} \frac{t^{\frac1q}}{S}\,.
\end{equation}
Similarly,
\begin{equation}
\label{tvt}
\frac{t'_v}{t} = q\cdot p^{\frac1q} \frac{t^{\frac1p}}{S}\,.
\end{equation}

Recall also that we had
\begin{equation}
\label{A}
A=A(u,v)=\frac{p^{\frac1p}}{p} t^{\frac1q}u + \frac{p^{\frac1q}}{q} t^{\frac1p}v\,.
\end{equation}

Using the notations \eqref{skobkaS} and \eqref{A} we can compute the Hessian of $B=B_q$. Namely,
$$
B_{uu} = \frac{2p^{\frac1p}}{pq} t^{\frac1q} \frac{t'_u}{t} -\frac1{pq} A (\frac{t'_u}{t})^2+\frac1{pq} S \frac{t''_{uu}}{t}\,.
$$
$$
B_{vv} = \frac{2p^{\frac1q}}{pq} t^{\frac1p} \frac{t'_v}{t} -\frac1{pq} A (\frac{t'_v}{t})^2+
+\frac1{pq} S \frac{t''_{vv}}{t}\,.
$$
$$
B_{uv} = \frac{p\,t}{S}-\frac1{pq} A \frac{t'_ut'_v}{t^2}+
+\frac1{pq} S \frac{t''_{uv}}{t}
$$
Plugging
$$
\frac{t''_{uu}}{t}= (\frac1q+1-\frac1p) (\frac{t'_u}{t})^2 -\frac{t}{S} (\frac{t'_u}{t})^2
$$ and using \eqref{tut} we get the following concise formulas:
\begin{equation}
\label{Buu1}
B_{uu} = \frac1{pq} S (\frac{t'_u}{t})^2\,.
\end{equation}

\begin{equation}
\label{Bvv1}
B_{vv} = \frac1{pq} S (\frac{t'_v}{t})^2\,.
\end{equation}

\begin{equation}
\label{Buu1}
B_{uv}+1  = \frac1{pq} S \frac{t'_ut'_v}{t^2}\,.
\end{equation}

Let us introduce the notations:
$$
\alpha=\frac{t'_u}{t}\,,\, \beta=\frac{t'_v}{t}\,, \,m=\frac1{pq}S\,,\, \tau =\frac{\alpha}{\beta}\,.
$$
Then we saw in the previous sections that the Hessian quadratic form of $B$
$$
Q(dx_1, dx_2) = B_{uu}dx_1^2 + 2(B_{uv} +1) dx_1dx_2 + B_{vv}dx_2^2
$$
will have the form
\begin{equation}
\label{Q}
Q= \frac{\alpha}{\beta} dx_1^2 + \frac{\beta}{\alpha} dx_2^2 +  \frac{\alpha}{\beta} (m\alpha\beta-1)(dx_1+\frac{\beta}{\alpha} dx_2)^2\,.
\end{equation}

\medskip

It is useful if the reader thinks that in what follows $y_{11}, y_{12}, y_{21}, y_{22}$ are, correspondingly, $\Phi, \Psi, U, V$.

Also in what follows $dy_{11}, dy_{12}, dy_{21}, dy_{22}$ can be viewed as $\phi_1, \psi_1, u_1, v_1$ and $\phi_2, \psi_2, u_2, v_2$.

\bigskip

Our goal now is to ``tensorize" the form $Q$. This operation means in our particular case to consider
the new function, now of $4$ real variables (or $2$ complex variables if one prefers), given by
$$
\B :=\B(y_{11}, y_{12}, y_{21}, y_{22}) := B(x_1,x_2)\,,\,\text{where}\,\, x_1:=\sqrt{y_{11}^2+y_{12}^{2}}\,,x_2:=\sqrt{y_{21}^2+y_{22}^{2}}
$$
and to write its Hessian quadratic form. In the previous section we saw the formula for doing that:
$$
\Q=  \frac{\alpha}{\beta}\bigg(\frac{y_{11}dy_{11} +y_{12}dy_{12}}{x_1}\bigg)^2 + \frac{\beta}{\alpha}\bigg(\frac{y_{21}dy_{21} +y_{22}dy_{22}}{x_2}\bigg)^2 +
$$
$$
\frac{\alpha}{\beta}(m\alpha\beta-1) \bigg(\frac{y_{11}dy_{11} +y_{12}dy_{12}}{x_1}+ \frac{\beta}{\alpha}\frac{y_{21}dy_{21} +y_{22}dy_{22}}{x_2}\bigg)^2 +
$$
$$
\frac{B_u}{u}\bigg(\frac{y_{12}dy_{11} -y_{11}dy_{12}}{x_1}\bigg)^2 + \frac{B_v}{v}\bigg(\frac{y_{22}dy_{21} -y_{21}dy_{22}}{x_2}\bigg)^2 \,.
$$

To show that this quadratic form has an interesting self-improving property we are going to make some calculations. 
First of all notice that
\begin{equation}
\label{albe}
\tau = \frac{\alpha}{\beta} = \frac{p\cdot p^{\frac1p}\cdot t^{\frac1q}}{q\cdot p^{\frac1q}\cdot t^{\frac1p}}
\end{equation}

Now we start with combining \eqref{Bu} with \eqref{tut}
\begin{equation}
\label{Buab}
B_u =\frac{p^{\frac1p}}{p}t^{\frac1q} -v +\frac1{pq} S\frac{p\cdot p^{\frac1p}\cdot t^{\frac1q}}{S}=p^{\frac1p}t^{\frac1q} -v\,.
\end{equation}

 Let us see that
 \begin{equation}
 \label{uff}
 \frac{p}{q}\frac{p^{\frac1p}}{p^{\frac1q}}\frac{t^{\frac1q}}{t^{\frac1p}} = \frac{p^{\frac1p}t^{\frac1q}}{u} -\frac{v}{u}\,.
 \end{equation}
 This is the same as
 $$
 p\cdot p^{\frac1p} \cdot t^{\frac1q} u = qp t - q\cdot p^{\frac1q}\cdot t^{\frac1p} v\,.
 $$
 
 But the last claim is correct, it is just the implicit equation \eqref{timpl} fot $t$.  So \eqref{uff} is correct. So, combining \eqref{albe} and \eqref{Buab} we obtain
 \begin{equation}
 \label{Buab1}
 \frac{B_u}{u} =\frac{\alpha}{\beta}\,.
 \end{equation}
 
 We would expect that $\frac{B_v}{v} =\frac{\beta}{\alpha} =\frac1{\tau}$ by symmetry, but actually $\frac{B_v}{v} >\frac{\beta}{\alpha}$ for $p>2$ and this allows us to have an improved inequality for $\Q$. Let us see how.
 
 Using \eqref{timpl} we get
 $$
 \frac{B_v}{v} -\frac{\beta}{\alpha} = \frac{p^{\frac1q} t^{\frac1p}}{v}-\frac{u}{v} - \frac{q\cdot p^{\frac1q}\cdot t^{\frac1p}}{p\cdot p^{\frac1p}\cdot t^{\frac1q}}=
 $$
 $$
 \frac{p^2 t - p\cdot p^{\frac1p}\cdot t^{\frac1q} u - q\cdot p^{\frac1q}\cdot t^{\frac1p}v}{p\cdot p^{\frac1p}\cdot t^{\frac1q}v}=
 $$
 $$
 \frac{(p^2-pq) t}{p\cdot  p^{\frac1p}\cdot t^{\frac1q}v} = \frac{p-q}{p^{\frac1p}} \frac{t^{\frac1p}}{v} =
 $$
 $$
 \frac{p-q}{p^{\frac1p}}\frac{\frac{p^{\frac1p}}{q} t^{\frac1q} u +\frac{p^{\frac1q}}{p} t^{\frac1p} v}{t^{\frac1q}v}=
 $$
 $$ 
 \frac{(1-q/p)}{p^{\frac1p-\frac1q}}  t^{\frac1p-\frac1q}+\bigg(\frac{p}{q}-1\bigg)\,.
 $$
 In particular, using \eqref{albe}
 $$
 \frac{B_v}{v} -\frac{\beta}{\alpha} +\frac{2}{\tau} \ge \frac{(1-q/p)}{p^{\frac1p-\frac1q}}  t^{\frac1p-\frac1q}+2\frac{q\cdot p^{\frac1q}\cdot t^{\frac1p}}{p\cdot p^{\frac1p}\cdot t^{\frac1q}}=
 $$
 $$
 \frac{(1-q/p)}{p^{\frac1p-\frac1q}}  t^{\frac1p-\frac1q} + \frac{2q/p}{p^{\frac1p-\frac1q}}  t^{\frac1p-\frac1q}=
 $$
 $$
\frac{q}{p^{\frac1p-\frac1q}}  t^{\frac1p-\frac1q} =p\cdot\frac1{\tau}\,.
$$
 
This is what we need
\begin{equation}
\label{Bvba}
\frac{B_v}{v} -\frac{\beta}{\alpha}  +\frac2{\tau} = p\cdot \frac1{\tau}+ (p/q-1) \frac{u}{v}\ge p\cdot \frac1{\tau}\,.
\end{equation}

\medskip

Now let us take a look at $\Q$ and let us plugg \eqref{Buab1} and \eqref{Bvba} into it. Then
\begin{equation}
\label{key}
\Q \ge \tau (dy_{11}^2 + dy_{12}^2)  + \frac1{\tau} (dy_{21}^2 + dy_{22}^2) +(\frac{B_v}{v} -\frac{\beta}{\alpha})\bigg(\frac{y_{22}dy_{21}-y_{21}dy_{22}}{x_2}\bigg)^2\,.
\end{equation}

Now imagine that we apply this estimate to {\bf two} different collection of vectors
$(dy_{11}, dy_{12}, dy_{21}, dy_{22})$, $(dy'_{11}, dy'_{12}, dy'_{21}, dy'_{22})$. Moreover, suppose that
we have orthonormality condition
\begin{equation}
\label{orth}
dy_{21}\cdot dy_{22}+dy'_{21}\cdot dy'_{22}=0\,,dy_{21}^2 + (dy'_{21})^2=dy_{22}^2+ (dy'_{22})^2\,.
\end{equation}
Then we get from \eqref{key}, \eqref{orth}
$$
\Q(dy) +\Q(dy') \ge \tau (dy_{11}^2 + dy_{12}^2+(dy'_{11})^2 + (dy'_{12})^2) + 1/\tau (dy_{21}^2 + dy_{22}^2+(dy'_{21})^2 + (dy'_{22})^2) +
$$
$$
(\frac{B_v}{v} -\frac{\beta}{\alpha}) \frac{y_{22}^2 +y_{21}^2}{x_2^2} \frac{(dy_{21}^2 + dy_{22}^2+(dy'_{21})^2 + (dy'_{22})^2) }{2}\,.
$$
We denote $\xi_1^2:=dy_{11}^2 + dy_{12}^2+(dy'_{11})^2 + (dy'_{12})^2, \xi_2^2:= dy_{21}^2 + dy_{22}^2+(dy'_{21})^2 + (dy'_{22})^2$.
Using that $\frac{y_{22}^2 +y_{21}^2}{x_2^2}=1$ and \eqref{Bvba} we rewrite the RHS and get
$$
\Q(dy) +\Q(dy') \ge \tau\cdot \xi_1^2 + \frac12 (\frac{B_v}{v} -\frac{\beta}{\alpha} +\frac2{\tau}) \xi_2^2\ge
\tau \cdot\xi_1^2 + 1/\tau \cdot\frac{p}2 \xi_2^2 \ge
$$
\begin{equation}
\label{key1}
 2\sqrt{\frac{p}2}(dy_{11}^2 + dy_{12}^2+(dy'_{11})^2 + (dy'_{12})^2)^{\frac12} 
(dy_{21}^2 + dy_{22}^2+(dy'_{21})^2 + (dy'_{22})^2)^{\frac12}\,.
\end{equation}

So we won $\sqrt{2/p}=\sqrt{\frac{2(q-1)}{q}}$ in comparison with the usual Burkholder estimate, which would be $\le \frac1{q-1}$. So the estimate for the orthogonal martingale will be
$\le \sqrt{\frac{2(q-1)}{q}}\cdot \frac1{q-1} =\sqrt{\frac{2}{q(q-1)}}$.

And we get Theorem \ref{less2}.

\markboth{}{\sc \hfill \underline{References}\qquad}

\end{document}